\documentclass[12pt]{article}
\usepackage[latin1]{inputenc}
\usepackage[english]{babel}
\usepackage{amsmath,amsthm,amscd,amssymb}
\usepackage{pb-diagram,bbm}
\usepackage{anysize}
\usepackage{multicol}
\usepackage{float}
\usepackage{graphics,graphicx}
\usepackage{caption}
\usepackage{undertilde}
\usepackage{subcaption}
\usepackage{cancel}
\numberwithin{equation}{section}

\theoremstyle{definition}

\marginsize{3cm}{2cm}{2.5cm}{2.5cm}

\title{{\bf Functional SAR Model}}
\author{Wilmer Pineda R\'ios \& Ram\'on Giraldo}
\begin{document}
\maketitle
\begin{abstract}

\end{abstract}

\section{Introduction}
In the last decades, advances in computer technology, modern equipment collection and storage of data, as well as advances in the different fields of science have enabled researchers to collect and provide data of high resolution digitized representing complex objects such as curves, surfaces or any element that varies on a continuous (time, space, wavelength, etc.). This is the case of data collected by seismographs, data on nuclear explosions, data on temperature, precipitation, medical data (electroencephalograms, electrocardiograms), financial data, which can be considered curves. A particular case of such data is that of variables taking values into an infinite dimensional space, typically a space of functions defined of some set $T$ \cite{Ramsay}.

The standard statistical techniques for modeling functional data are focused on independent functions. However, in several disciplines of applied sciences there exists an increasing interest in modeling correlated functional data: this is the case when samples of functions are observed over a discrete set of time points (temporally correlated functional data) or when these functions are observed in different sites of a region (spatially correlated functional data). In these cases the above-mentioned methodologies may not be appropriate as they do not incorporate dependence among functions into the analysis. For this reason, some statistical methods for modeling correlated variables, such as time series (Box and Jenkins 1976) or geostatistical analysis \cite{Cressie}, have been adapted to the functional context.

The plan of the article is as follows. Section 2 presents the Simultaneous Autoregressive (SAR) Model in the case where the explanatory variable is continous. Section 3 introduces the functional SAR model. Section 4 gives a simulation study. The article ends with a brief discussion and suggestions for further research.

\section{Simultaneous Autoregressive (SAR) Model}
\noindent In the context of time series, autoregressive model represents the dependence of observation at time $t$ as a linear combination of its past values. A spatial analog may be defined in the following sense: to represent the dependence of observation in the area $\mathbf{s}$ as a linear combination of its neighboring observations. This consideration allows for spatial dependence in the data. If we apply this idea spatial autoregression on residual vector, the following structure is obtained:

\begin{align} \label{SAR}
 Z(\mathbf{s}) &= X(\mathbf{s})\beta + e(\mathbf{s})\\
 e(\mathbf{s}) &= Be(\mathbf{s}) + v
\end{align}

\noindent where $B$ is a matrix of spatial dependence parameters with $b_{ii}=0$. In the literature, it is considered that $v$ have mean zero and a diagonal covariance matrix $\Sigma$. If all $b_{ij}$ are zero, there is no autoregression and the model reduces to the traditional linear regression model with uncorrelated errors.

\noindent If we solve \eqref{SAR} for $v$, this autoregressive model can be express as

\begin{equation}\label{SAR1}
 (I-B)\left(Z(\mathbf{s})- X(\mathbf{s})\beta\right) = v
\end{equation}

\noindent The model in \eqref{SAR1} was introduced by Whittle in \cite{Whittle} and the adjective ``simultaneous'' describes the $n$ autoregression that occur simultaneously at each data location in this formulation. The matrix of spatial dependence parameters $B$ plays an important role in SAR models. In \cite{Schabenberger} propose to take $B=\rho W$, where $W$ is a known spatial proximity matrix, in order to make progress with estimation and inference. With this parametrization of $B$, the SAR model can be written as

\begin{align*}
 Z(\mathbf{s}) &= X(\mathbf{s})\beta + e(\mathbf{s})\\
  e(\mathbf{s}) &= \rho We(\mathbf{s}) + v
\end{align*}

\noindent For a well-defined model, it is necessary that $(I-\rho W)$ to be an invertible matrix. This restriction imposes conditions on $W$ and also on $\rho$. Haining in \cite{Haining} concluded that if $\lambda_{\max}$ and $\lambda_{\min}$ are the largest and smallest eigenvalues of $W$, and if $\lambda_{\min}<0$ and $\lambda_{\min}>0$, then
\begin{equation*}
 \frac{1}{\lambda_{\min}}<\rho<\frac{1}{\lambda_{\max}}
\end{equation*}

\noindent Often, the row sums of $W$ are standardized to 1 by dividing each entry in $W$ by its row sum, Then $\lambda_{\max}=1$ and $\lambda_{\min}\leq -1$, so $\rho<1$. If $\rho$ is known and $\Sigma=\sigma^2I$ then generalized least squares can be used to estimate $\beta$ and $\sigma^2$. Thus,

\begin{align*}
 \hat{\beta} &= \left(X(\mathbf{s})\Sigma_Z^{-1}X(\mathbf{s})\right)^{-1}(X(\mathbf{s})\Sigma_Z^{-1}Z(\mathbf{s})\\
 \hat{\sigma}^2 &= \frac{\left(Z(\mathbf{s})-X(\mathbf{s})\hat{\beta}\right)^T\Sigma_Z^{-1}\left(Z(\mathbf{s})-X(\mathbf{s})\hat{\beta}\right)}{n-k}\\
 \Sigma_Z &= \sigma^2(I-\rho W)^{-1}(I-\rho W^T)^{-1}
\end{align*}

\section{Functional SAR Model}
\noindent In this section, a functional SAR model with scalar response is proposed. The data we observe for the $i$th area are $\{\left(X_i(t); t\in\mathcal{T},Y_i\right)\}$. The predictor variable $X(t)$, $t\in\mathcal{T}$, is a random curve which is observed per area of experimental unity and corresponds to a square integrable stochastic process on a real interval $\mathcal{T}$. The dependent variable $Y$ is a real-valued continous random variable. In this model, a structure of spatial autocorrelation is assigned to the residual of the scalar response model introduced in \cite{Kokoszka}, using a matrix of spatial proximity between areas. The functional SAR model is proposed as
\begin{equation} \label{E1}
 \begin{cases}
  Y=\int\limits_{\mathcal{T}} \mathbf{X}(t)\beta(t) dt + \nu\\
  \nu = \rho W \nu + \varepsilon
 \end{cases}
\end{equation} 

\noindent where $\varepsilon\sim N(0,\sigma^2I)$ and $W$ is a symmetric proximity matrix. We express the model \eqref{E1} as follow

\begin{align*}
 Y &= \int\limits_{\mathcal{T}} \mathbf{X}(t)\beta(t) dt + \rho W \nu + \varepsilon\\
 &= \int\limits_{\mathcal{T}} \mathbf{X}(t)\beta(t) dt + \rho W\left(Y - \int\limits_{\mathcal{T}} \mathbf{X}(t)\beta(t) dt\right) + \varepsilon\\
 &= \int\limits_{\mathcal{T}} \mathbf{X}(t)\beta(t) dt + \rho WY - \rho W\int\limits_{\mathcal{T}} \mathbf{X}(t)\beta(t) dt + \varepsilon
\end{align*}

Solving for $\varepsilon$,

\begin{equation}\label{Equation2}
 \varepsilon = (I-\rho W)\left(Y - \int\limits_{\mathcal{T}} \mathbf{X}(t)\beta(t) dt\right)
\end{equation}

The parameter function $\beta(t)$ is a quantity of central interest in the statistical analysis and replaces the vector of slopes in a linear model. Let $\{\varphi_j\}_{j\geq 1}$ be an orthonormal basis of the function space $L^2(\mathcal{T})$. Then the predictor process $X(t)$ and the parameter function $\beta(t)$ can be expanded into

\begin{equation*}
 X_i(t)=\sum\limits_{j=1}^{\infty} a_{ij}\varphi_j(t), \qquad \beta(t)=\sum\limits_{j=1}^{\infty} b_j\varphi_j(t)
\end{equation*}

Here, $\{a_{ij}\}_{j\geq 1}$ are random variables associated with the $i$th area. We note that, using Parseval identity,

\begin{align*}
 \sum\limits_{j=1}^{\infty} b_j^2 &= \sum\limits_{j=1}^{\infty}\left|\langle \beta(t), \varphi_j(t)\rangle\right|_{L^2}^2\\
 &= \|\beta(t)\|_{L^2}^2 < \infty 
\end{align*}

Furthermore, if $\mathbb{E}(a_{ij})=0$ and $\mathbb{E}\left(a_{ij}^2\right)=\sigma_j^2$, for all $i=1,\dotsc,n$,

\begin{equation*}
\sum\limits_{j=1}^{\infty} \sigma_j^2 = \int\limits_{\mathcal{T}} \mathbb{E}\left(X_i^2(t)\right) dt < \infty
\end{equation*}

Now,

\begin{align*}
 \int\limits_{\mathcal{T}} X_i(t)\beta(t) dt &= \int\limits_{\mathcal{T}} \left(\sum\limits_{k=1}^{\infty} a_{ik}\varphi_k(t)\right) \left(\sum\limits_{j=1}^{\infty} b_j\varphi_j(t)\right) dt\\
 &= \sum\limits_{k=1}^{\infty}\sum\limits_{j=1}^{\infty} a_{ik}b_j\left(\int\limits_{\mathcal{T}} \varphi_k(t)\varphi_j(t) dt\right)\\
 &=\sum\limits_{j=1}^{\infty} a_{ij}b_j
\end{align*}

The model proposed here has a difficulty caused by the infinite \-di\-men\-sio\-na\-li\-ty of $L^2(\mathcal{T})$. So, the model is truncated at $k=k_n$ and the dimension $k_n$ increases asymptotically as $n\to\infty$. A truncation strategy is as follow:

\begin{equation*}
 \begin{cases}
  U_{ik}=\sum\limits_{j=1}^{k} a_{ij}b_j\\
  V_{ik}=\sum\limits_{j=k+1}^{\infty} a_{ij}b_j
 \end{cases}
\end{equation*}

Let, $\mathbf{A}=(a_{ij})$ is a $k\times k$ matrix and $\mathbf{b}_k=(b_1,\dotsc,b_k)^T$ is a $k$-dimensional vector of slopes. If $U=\mathbf{A}\mathbf{b}_k$, Let $g(U)=\mathbb{E}(Y|U)$ be the conditional expected value and let $F_{V|U}$ be the conditional distribution of $V=\int\limits_{\mathcal{T}} \mathbf{X}(t)\beta(t) dt - U$ given $U$. Hence,

\begin{align*}
 \mathbb{E}\left(\left(U+V-g(U)\right)^2\right) &= \mathbb{E}\left(\left(\int (U+V)-(U+s) dF_{V|U}\right)^2\right)\\
 &= \mathbb{E}\left(\left(\int (V-s) dF_{V|U}\right)^2\right)\\
 &\leq 2\mathbb{E}\left(\int (V^2+s^2) dF_{V|U}\right)\\
 &= 2\mathbb{E}\left(\mathbb{E}\left(V^2|U\right)+\mathbb{E}\left(V^2\right)\right)\\
 &= 4\mathbb{E}\left(V^2\right)\\
 &= 4\left(\sum\limits_{j=k+1}^{\infty} b_j^2\right)\left(\sum\limits_{j=k+1}^{\infty} \sigma_j^2\right)
\end{align*}

Then, the approximation error of truncated model is seen to be directly tied  to $\mathbf{Var}(V)$ and is controlled by the sequence $\sigma_j^2=\mathbf{Var}(a_{ij})$, $j=1,2,\dotsc$, corresponding to a sequence of eigenvalues for the special case of a eigenbase. Inference will be developed using asymptotic results ($k\to \infty$).

\subsection{Least Squares Estimation}

One central goal is estimation and inference for $\beta(t)$. In this case, using the truncated model, we define

\begin{equation*}
 \varepsilon_k= (I-\rho_k W)\left(Y - \mathbf{A}\mathbf{b}_k\right)
\end{equation*}

The idea is finding $\mathbf{b}_k$ and $\rho_k$ that minimizes the following expression

\begin{align*}
 \varepsilon_k^T\varepsilon_k &= \left(Y - \mathbf{A}\mathbf{b}_k\right)^T(I-\rho_k W)^2\left(Y - \mathbf{A}\mathbf{b}_k\right)\\
 &= \left(Y - \mathbf{A}\mathbf{b}_k\right)^TZ\left(Y - \mathbf{A}\mathbf{b}_k\right)
\end{align*}

where $Z=(I-\rho_k W)^2$. Hence,

\begin{align*}
 \varepsilon_k^T\varepsilon_k &= Y^TZY - Y^TZ\mathbf{A}\utilde{\mathbf{b}_k} - \mathbf{b}_k^T\mathbf{A}^TZY + \mathbf{b}_k^T\mathbf{A}^TZ\mathbf{A}\mathbf{b}_k\\
 &= Y^TZY - 2\mathbf{b}_k^T\mathbf{A}^TZY + \mathbf{b}_k^T\mathbf{A}^TZ\mathbf{A}\mathbf{b}_k
\end{align*}

Taking derivative of the last equality with respect to $\mathbf{b}_k$,

\begin{equation} 
\frac{d\varepsilon_k^T\varepsilon_k}{d\mathbf{b}_k} = -2\mathbf{A}^TZY+2\mathbf{A}^TZ\mathbf{A}\mathbf{b}_k
\end{equation}

Now, taking derivative of $\varepsilon_k^T\varepsilon_k$ with respect to $\rho_k$,

\begin{align*}
\frac{d\varepsilon_k^T\varepsilon_k}{d\rho_k} &= Y^T\left(2\rho_kW^2-2W\right)Y - 2\mathbf{b}_k^T\mathbf{A}^T\left(2\rho_kW^2-2W\right)Y + \mathbf{b}_k^T\mathbf{A}^T\left(2\rho_kW^2-2W\right)\mathbf{A}\mathbf{b}_k\\
&= 2\rho_k\left(Y^TW^2Y-2\mathbf{b}_k^T\mathbf{A}^TW^2Y+\mathbf{b}_k^T\mathbf{A}^TW^2\mathbf{A}\mathbf{b}_k\right) - 2\left(Y^TWY-2\mathbf{b}_k^T\mathbf{A}^TWY+\mathbf{b}_k^T\mathbf{A}^TW\mathbf{A}\mathbf{b}_k\right)\\
&= 2\rho_k\left(Y-\mathbf{A}\mathbf{b}_k\right)^TW^2\left(Y-\mathbf{A}\mathbf{b}_k\right) - 2\left(Y-\mathbf{A}\mathbf{b}_k\right)^TW\left(Y-\mathbf{A}\mathbf{b}_k\right)
\end{align*}

Hence, we obtain the following system of equations

\begin{equation}
 \begin{cases}
  \mathbf{A}^T\left(I-\hat{\rho}_kW\right)^2\mathbf{A}\mathbf{\hat{b}}_k = \mathbf{A}^T\left(I-\hat{\rho}_kW\right)^2Y \\
  \rho_k\left(Y-\mathbf{A}\mathbf{\hat{b}}_k\right)^TW^2\left(Y-\mathbf{A}\mathbf{\hat{b}}_k\right) = \left(Y-\mathbf{A}\mathbf{\hat{b}}_k\right)^TW\left(Y-\mathbf{A}\mathbf{\hat{b}}_k\right)
 \end{cases}
\end{equation}

But the above system cannot be solved analytically. Then, for finding $\hat{\rho}_k$ and $\mathbf{\hat{b}}_k$, we establish an iterative procedure according to \cite{Schabenberger}:

\begin{enumerate}
\item Estimate $\mathbf{b}_k$, taking $\rho_k=0$
\item Estimate $\rho_k$ using the last estimator for $\mathbf{b}_k$
\item Estimate $\mathbf{b}_k$ using the last estimator for $\rho_k$
\item Repeat steps 2 and 3 until convergence.
\end{enumerate}

\subsection{Maximum Likelihood Estimation}

Now, we want to find the maximum likelihood estimator for $\beta(t)$ and $\rho$. Using \eqref{Equation2} and the assumption over $\varepsilon$,

\begin{equation*}
 \mathbb{E}\left((I-\rho W)\left(Y-\int\limits_{\mathcal{T}}X(t)\beta(t) dt\right)\right) =0
\end{equation*}

If $(I-\rho W)$ is an invertible matrix, using properties of expected value,

\begin{equation*}
 \mathbb{E}(Y)= \int\limits_{\mathcal{T}}X(t)\beta(t) dt
\end{equation*}

For the covariance matrix of $Y$,

\begin{align*}
 Cov(Y) &= Cov(v)\\
 &= Cov\left((I-\rho W)^{-1}\varepsilon\right)\\
 &= (I-\rho W)^{-1}Cov(\varepsilon)(I-\rho W)^{-1}\\
 &=\sigma^2\left((I-\rho W)^2\right)^{-1}
\end{align*}

Hence,

\begin{equation*}
 Y\sim N\left(\int\limits_{\mathcal{T}}X(t)\beta(t) dt, \sigma^2\left((I-\rho W)^2\right)^{-1}\right)
\end{equation*}

The likelihood function is given by

\begin{align*}
 L(\beta(t),\sigma^2,\rho) &= (2\pi)^{-n/2}(\sigma^2)^{-n/2}\left|(I-\rho W)^2\right|^{1/2}\\ & \exp\left\{-\frac{1}{2\sigma^2}\left(Y-\int\limits_{\mathcal{T}}X(t)\beta(t) dt\right)^T(I-\rho W)^2\left(Y-\int\limits_{\mathcal{T}}X(t)\beta(t) dt\right)\right\}
\end{align*}

Taking natural logarithm, we obtain (without the constant) the function to maximize

\begin{align*}
 l(\beta(t),\sigma^2,\rho) &= -\frac{n}{2}\ln(\sigma^2) + \ln(\left|(I-\rho W)\right|)\\  &- \frac{1}{2\sigma^2}\left(Y-\int\limits_{\mathcal{T}}X(t)\beta(t) dt\right)^T(I-\rho W)^2\left(Y-\int\limits_{\mathcal{T}}X(t)\beta(t) dt\right)\\
 &= -\frac{n}{2}\sigma^2 + \ln(\left|(I-\rho W)\right|) - \frac{1}{2\sigma^2}\left(Y-\mathbf{A}\mathbf{b}_k-V\right)^T(I-\rho W)^2\left(Y-\mathbf{A}\mathbf{b}_k-V\right)\\
\end{align*}

Again, using the truncated model, the function to maximize is approximate to

\begin{equation*}
 l_k(\beta(t),\sigma^2,\rho) = -\frac{n}{2}\ln(\sigma^2) + \ln(\left|(I-\rho_k W)\right|) - \frac{1}{2\sigma^2}\left(Y-\mathbf{A}\mathbf{b}_k\right)^T(I-\rho_k W)^2\left(Y-\mathbf{A}\mathbf{b}_k\right)
\end{equation*}

Then,

\begin{align*}
  \frac{\partial l_k}{\partial\mathbf{b}_k} &= -\frac{1}{2\sigma^2}\frac{d\varepsilon_k^T\varepsilon_k}{d\mathbf{b}_k}\\
  \frac{\partial l_k}{\partial\rho_k} &= -Tr\left((I-\rho_k W)^{-1}W\right) - \frac{d\varepsilon_k^T\varepsilon_k}{d\rho_k}\\
  \frac{\partial l_k}{\partial\sigma^2} &= -\frac{n}{2\sigma^2} + \frac{1}{2\sigma^4}\left(Y-\mathbf{A}\mathbf{b}_k\right)^T(I-\rho_k W)^2\left(Y-\mathbf{A}\mathbf{b}_k\right)
\end{align*}

Hence, the log-likelihood function should be maximized using numerical methods. Using any method of estimation, we write,

\begin{equation}
 \hat{\beta}_k(t)=\mathbf{\Phi}_k(t)^T\mathbf{\hat{b}}_k
\end{equation}

where $\mathbf{\Phi}_k(t)=\left(\varphi_1(t),\dotsc,\varphi_k(t)\right)^T$ is the $k$-dimensional vector of the first $k$ elements of the initial basis of $L^2(\mathcal{T})$.

\subsection{Confidence Bands and Hyphothesis Testing}
Now, let us assume that $\rho$ is known and $k=k_n$ is fixed, in order to make inference (confidence bands and hypotheses) about $\beta(t)$. We write

\begin{equation*}
 \beta_k(t)=\mathbf{\Phi}_k(t)^T\mathbf{b}_k
\end{equation*}

And let $\Sigma=\mathbf{A}^T(I-\rho W)^2\mathbf{A}$. Notice that $\Sigma$ is a symmetric matrix with real entries, hence $\Sigma$ can be diagonalized. Then,

\begin{equation*}
 \Sigma=PJP^{-1}
\end{equation*}
where $J$ is the diagonal matrix of eigenvalues of $\Sigma$. Using the assumption ($\rho$ known), 

\begin{equation*}
 \mathbf{\hat{b}}_k = \Sigma^{-1}\mathbf{A}^T\left(I-\rho W\right)^2Y
\end{equation*}

Then, computing the expected conditional value of $\mathbf{\hat{b}}_k$ given $U$,

\begin{align*}
 \mathbb{E}\left(\mathbf{\hat{b}}_k|U\right) &= \Sigma^{-1}\mathbf{A}^T\left(I-\rho W\right)^2\mathbb{E}\left(Y|U\right)\\
 &= \Sigma^{-1}\mathbf{A}^T\left(I-\rho W\right)^2\mathbf{A}\mathbf{b}_k\\
 &= \Sigma^{-1}\Sigma\mathbf{b}_k\\
 &= \mathbf{b}_k
\end{align*}

Hence, $\mathbf{\hat{b}}_k$ is unbiased estimator for $\mathbf{b}_k$. On the other hand,

\begin{align*}
 Var\left(\mathbf{\hat{b}}_k|U\right) &= \Sigma^{-1}\mathbf{A}^T(I-\rho W)^2 Var\left(Y|U\right)(I-\rho W)^2\mathbf{A}\Sigma^{-1}\\
 &= \sigma^2\Sigma^{-1}\mathbf{A}^T(I-\rho W)^2 (I-\rho W)^{-2}(I-\rho W)^2\mathbf{A}\Sigma^{-1}\\
 &= \sigma^2\Sigma^{-1}\mathbf{A}^T(I-\rho W)^2\mathbf{A}\Sigma^{-1}\\
 &= \sigma^2\Sigma^{-1}\Sigma\Sigma^{-1}\\
 &= \sigma^2\Sigma^{-1}
\end{align*}

Using the results in \cite{Schabenberger}, the generalized least squares can be used to estimate $\sigma^2$. Thus,

\begin{equation}
 \hat{\sigma}^2 = \frac{\left(Y-\mathbf{A}\mathbf{\hat{b}}_k\right)^T(I-\rho W)^2\left(Y-\mathbf{A}\mathbf{\hat{b}}_k\right)}{n-r}
\end{equation}

where $r=\text{trace}\left(\mathbf{A}\Sigma^{-1}\mathbf{A}^T\left(I-\rho W\right)^2\right)$. Hence, a confidence band of \\ $100(1-\alpha)\%$ is

\begin{equation*}
 \hat{\beta}_k(t) \pm z_{1-\frac{\alpha}{2}}\hat{\sigma}\sqrt{\mathbf{\Phi}_k(t)^T\Sigma^{-1}\mathbf{\Phi}_k(t)}
\end{equation*}

Now, an important problem in practical situations is testing the following hypothesis

\begin{equation}\label{H1}
 \begin{cases}
  H_0: \beta(t)=\beta_0(t)\\
  H_a: \beta(t)\neq\beta_0(t)
 \end{cases}
\end{equation}

where $\beta_0(t)$ is a known function in $L^2(\mathcal{T})$. \cite{Cardot02} introduce two test statistics based on the square norm of a normalized version of the cross-covariance. We use and extend these results for the model introduced here. The model \eqref{E1} can be express as follow

\begin{equation}\label{E2}
 Q_i = \int\limits_{\mathcal{T}} Z_i(t)\beta(t)dt + \varepsilon_i
\end{equation}

where $Q=(I-\rho W)^{-1}Y$ is a modification of the response variable and

\begin{equation*}
 Z_i(t) = \sum\limits_{j=1}^{n} c_{ij}X_j(t)
\end{equation*}

Here, $c_{ij}$ is the $i-j$ element of the matrix $(I-\rho W)^{-1}$. The empirical covariance and cross operators for the model \eqref{E2} are defined by

\begin{align*}
 \Gamma_n x(t) &= \frac{1}{n}\sum\limits_{i=1}^{n} \langle Z_i(t),x(t)\rangle Z_i(t)\\
 \Delta_n x(t) &= \frac{1}{n}\sum\limits_{i=1}^{n} \langle Z_i(t),x(t)\rangle Q_i
\end{align*}

and let $(\hat{\lambda}_j,\hat{V}_j)$ be the eigenelements of $\Gamma_n$. Rewitring the expression for $\Delta_n$,

\begin{align*}
 \Delta_n x(t) &= \frac{1}{n}\sum\limits_{i=1}^{n} \langle Z_i(t),x(t)\rangle Q_i\\
 &= \frac{1}{n}\sum\limits_{i=1}^{n}\sum\limits_{j=1}^{n} c_{ij} \langle X_j(t),x(t)\rangle Q_i\\
 &= \frac{1}{n}\sum\limits_{i=1}^{n}\sum\limits_{j=1}^{n} c_{ji} \langle X_j(t),x(t)\rangle \left(\sum\limits_{k=1}^{n} c_{ik} Y_k\right)\\
 &= \frac{1}{n}\sum\limits_{j=1}^{n}\sum\limits_{k=1}^{n} \left(\sum\limits_{i=1}^{n} c_{ji}c_{ik}\right) \langle X_k(t), x(t)\rangle Y_j
\end{align*}

Using the notation used in \cite{Cardot02}, let

\begin{equation*}
 \hat{A}_n(\cdot) = \sum\limits_{j=1}^{k_n} \hat{\lambda}_j^{-\frac{1}{2}}\langle \hat{V}_j,\cdot\rangle\hat{V}_j
\end{equation*}

and

\begin{equation}
 T_n = \frac{1}{\sqrt{k_n}}\left(\frac{1}{\hat{\sigma}^2}\|\sqrt{n}\Delta_n\hat{A}_n\|^2 - k_n\right)
\end{equation}

the test statistic for the hypothesis testing \eqref{E1} with $\beta_0(t)=0$. $T_n$ converges in distribution, under $H_0$ to a centred gaussian real random variable with variance 2. Hence, given a significance level $\alpha$, let $z_{\frac{\alpha}{2}}$ be the quantile of order $\frac{\alpha}{2}$ of a normal distribution with mean 0 and variance 1. The rejection criterion of $H_0$ is $|T_n|>\sqrt{2}z_{\frac{\alpha}{2}}$.\\

In \cite{Cardot02}, they clarify that for the general case of testing $\beta(t)=\beta_0(t)$, only need to center the response the response variable of the following form

\begin{equation*}
 \tilde{Q}=Q-\langle \beta_0(t),Z\rangle
\end{equation*}

\section{A simulation study}

\begin{figure}[H]
\centering
\includegraphics[scale=0.5]{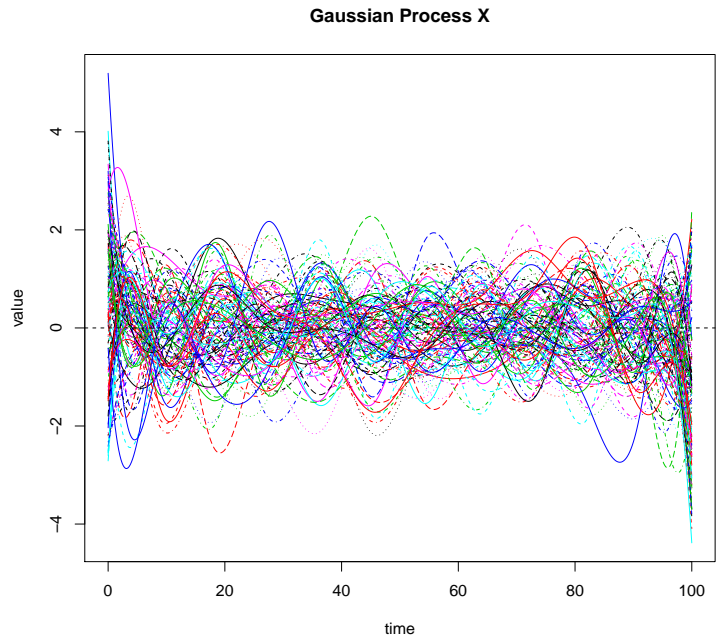}
\caption{\emph{Simulation of the Gaussian process $X(t)$ with systematic sinusoidal trend for 117 municipalities of Cundinamarca, Colombia}}
\label{fig1}
\end{figure}

\begin{figure}[H]
\centering
\includegraphics[scale=0.5]{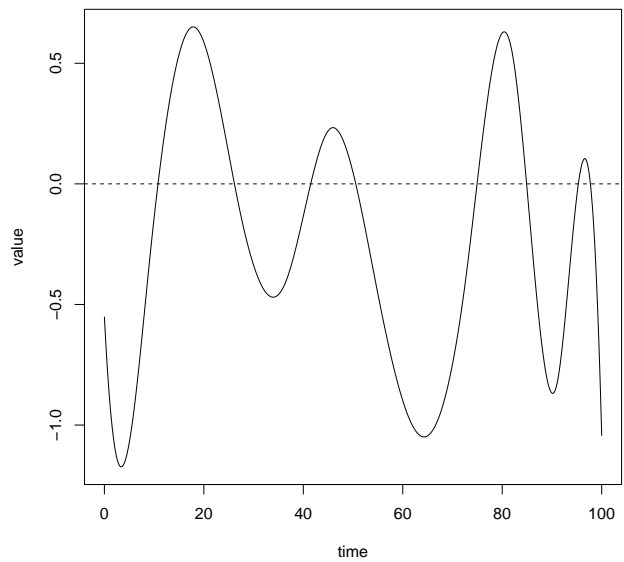}
\caption{\emph{Simulation of the parameter function $\beta(t)$}}
\label{fig2}
\end{figure}

We have considered five scenarios to assess the performance of the algorithm proposed. These scenarios are computed from a Gaussian process evaluated in a fine grid of $N=101$ points $\{t_1,\dotsc, t_n\}\in [0,100]$. Also, we added a systematic sinusoidal trend. The spatial weight matrix was built based on the coordinates of the departament of Cundinamarca, Colombia and $\rho$ and the spatial dependence parameter $\rho$ was taken as $0.1$, $0.5$ and $0.9$. Additionally, $\beta(t_i)=\cos(2t_i)+\varepsilon_i$, where $\varepsilon_i\sim N(0,2)$. For smoothing curves $X_i(t)$, $i=1,\dotsc, n$ and $\beta(t)$, it was used a base of cubic splines and the procedure was done using the {\it fda} package in R. In the figures \ref{fig1} y \ref{fig2} we show the results of the simulations for the Gaussian process $X(t)$ and $\beta(t)$.

Now, we estimate the functional SAR model with the above considerations. In the table \ref{tab1} we show the results for the five different scenarios for $\rho$, using $m=100$ distinct values of the vector Y and take the mean of the $\hat{\rho}$, the estimaton of expected value of integrated mean square error ($MISE$) associated to $\hat{\beta}(t)$ and the estimation of the expected value of $\sigma^2$.

\begin{table}[H]
\centering
\caption{\emph{Summary statistics of estimations for $\rho$, $MISE$ and $\sigma^2$}}
    \begin{tabular}{|c|c|c|c|}
    \hline
    $\rho$   & $\hat{\rho}$ & $MISE$  & $\hat{\sigma}^2$ \\
    \hline
    0.1   & 0.091 & 0.131 & 0.872 \\
    0.3   & 0.305 & 0.134 & 0.894 \\
    0.5   & 0.527 & 0.117 & 0.871 \\
    0.7   & 0.709 & 0.119 & 0.841 \\
    0.9   & 0.908 & 0.114 & 0.865 \\
    \hline
    \end{tabular}%
\label{tab1}
\end{table}

At least in terms of this particular type of data, the simulation suggests that the implemented algorithm detects in a good grade of accuracy the spatial dependence parameter $\rho$, but $\sigma^2$ is underestimated, which allows us to conclude that an approach may need a larger base of B-splines. The integrated mean square error is approximately 12\% and is expected to decrease when the functional variable is approximated with more elements in the B-Splines base and, therefore, the parameter function $\beta(t)$.

\end{document}